\documentstyle{article}
\newcommand{\bewende}{\hspace*{\fill}$\nabla $\\ \vspace{1cm}} 
\newcommand{\bewanf}{\noindent {\bf Proof:}\quad }    
\newcommand{\qua}{\enspace}                            
\newcommand{\avec}[2]{\mbox{$#1_1\ldots #1_#2$}}     

\newcommand{\prae}[2]{\mbox{$<#1\qua |\qua #2>$}}

\newcommand{\maN}{\mathop{\rm I\! N}\nolimits}

\newcommand{\maR}{\mathop{\rm I\! R}\nolimits}
\newcommand{\satzanf}{\begin{samepage}\begin{satz}}
\newcommand{\satzende}{\end{satz}\end{samepage}}
\newcommand{\koranf}{\begin{samepage}\begin{kor}}
\newcommand{\korende}{\end{kor}\end{samepage}}
\newcommand{\lemanf}{\begin{samepage}\begin{lem}}
\newcommand{\lemende}{\end{lem}\end{samepage}}
\newcommand{\bspanf}{\begin{beispiel}}
\newcommand{\bspende}{\end{beispiel}}
\newcommand{\propanf}{\begin{samepage}\begin{prop}}
\newcommand{\propende}{\end{prop}\end{samepage}}
\newtheorem{satz}{Theorem}[section]             
\newtheorem{lem}[satz]{Lemma}
\newtheorem{beispiel}[satz]{Example}
\newtheorem{kor}[satz]{Corollary}
\newtheorem{prop}[satz]{Proposition}
\newcommand{\be}{\begin{equation}}
\newcommand{\ee}{\end{equation}}
\newcommand{\bl}[1]{\begin{equation}\label{#1}}
\newcommand{\cay}{$\Gamma_X(G)$}
\newcommand{\cayq}{$\Gamma_X(G_q)$}
\newcommand{\normclz}{<\! <z>\! >}
\newcommand{\dist}{$\delta (\sigma ,w,v)$}
\begin{document}
\title{A bicombing that implies a\\
       sub -- exponential Isoperimetric Inequality}
\author{
        G\"unther Huck\\
        Northern Arizona\\University
     \and
        Stephan Rosebrock\\
        Johann--Wolfgang Goethe\\Universit\"at}

\maketitle

\section{Definitions}

    Let \cay\ be the Cayley graph
of a group $G$ with respect to a finite set of generators $X$,
and let \cay\ be equipped with the word metric.
Let $F$ be the free group on $X$. For $v\in F$ let $|v|$ denote the length
in the free group.

A bicombing as defined in \cite{AB} and \cite{Sh1} is essentially a
selection of a path
$\sigma(g,h)$ for every pair of vertices $g,h\in $\cay , such that the
distance between any two paths which start and end a distance $\le 1$
apart is
uniformly bounded. We replace the uniform bound for this distance by a bound
that is dependent on the lengths of the paths. More precisely, we define a
bicombing of narrow shape as follows:

   For each $(g,h) \in G \times G$ let $\sigma(g,h): [0,\infty [ \to$\cay\ be a path from
$g$ to $h$ which is at integer times at vertices (i.e.~from $t=n$ to
$t=n+1$ the path
either travels the distance between two adjacent vertices or pauses at a vertex).
We define the length:\\
$|\sigma(g,h)| = \min \{t|\sigma(g,h)[t,\infty [\qua = \mbox{constant} =h\}$.
This is the length of the path including the pauses which occur before the
end of the path is reached. We will frequently represent such
a path by a sequence of elements in $X \cup X^{-1} \cup \{1\}$ which, given the
startvertex $g$, completely determines the path.
Let $\sigma(h) = \sigma(1,h)$.
We call $\sigma $ a {\it bicombing of narrow shape} if
\begin{enumerate}
\item it is ``recursive'', i.e.~if
there exists an increasing polynomial $f\colon \maN \to \maN$, such that
\begin{equation}\label{konsm}  |\sigma(g)|\le f(d(1,g)) \qquad
\forall g\in G   \end{equation}

\item there exists an integer $M>1$ and a real number $k>2$, such that
for any \mbox{$g\in G$} $|\sigma (g,g)|\le Mk/2$ and
for all $g,h\in G$ and $a,b\in X^{\pm 1}\cup \{ 1\} $
\begin{equation}\label{konsk} |\sigma (\sigma (g,h)(t),\sigma (ga,hb)(t))|
\le \max((|\sigma (g,h)|+|\sigma (ga,hb)|)/k,M/2) \end{equation}
holds for all integers $t\in [0,\infty [ $.
\end{enumerate}
where $d(1,g)$ denotes the distance in \cay\ from 1 to $g$. If possible
we will always choose $\sigma (1)$ to be the identical path.
A bicombing is called {\it geodesic} if $f$ is the identity
(i.e.~the combing lines are geodesics).

    Let the group $G$ be finitely generated with generator set $X$.
Following Gersten \cite{G8}, a function $f: \maN \to \maN$
is called an {\it isoperimetric function} for
$G$ if for any word $w$ in $X$ of length $n$ with $w=1$ in $G$, the
minimum number of 2-cells in a van Kampen diagram for $w$ is at
most $f(n)$.

    Let $P=\prae{X}{R}$ be a finite presentation of the group $G$.
Following Gersten \cite{G8}, a function $f\colon \maN \to \maR$ is
called an {\it isodiametric function} for $P$, if for any word $w$ in the
generators $X$ with $w=1$ in $G$ there is a van Kampen diagram for $w$,
such that any vertex in the diagram has distance at most
$f(|w|)$ from the basepoint.

   We would like to thank Allan J.~Sieradsky, Holger Meinert, Stephen J.~Pride, William
A.~Bogley and all the members of our ''Luttach workshop'' for
helpful discussions.

\section{An isoperimetric inequality and an isodiametric function}

\satzanf \label{satzIsop} A group $G$ with finite generator set $X$ and a
bicombing of narrow shape is finitely presented and has an
isoperimetric function of growth $n^{O(\log n)}$. \satzende

\bewanf Define a presentation $P=\prae{X}{R}$, where $R$ is the set of
all cyclically reduced non-trivial words of length at most $M+2$
which are trivial in $G$. We proof that $P$ is a presentation for
$G$ by constructing a van Kampen diagram for each
word which is trivial in $G$, using only 2--cells of $R$.

      Let $w\in F$ be a reduced nontrivial word of length $n>M+2$
which is trivial in $G$. If $w=\avec{x}{n},\quad x_i\in X^{\pm 1}$,
define $w_i=\avec{x}{i}$. Now consider the
''fan`` of bicombing lines $\sigma(w_i)$ from 1 to $w_i$. $|w|=n$ implies
$d(1,w_i)\le n/2$ and by (\ref{konsm}) it follows
\be\label{konsnm} |\sigma (w_i)|\le f(n/2) \quad \mbox{for }1\le i\le n.\ee

    If $|\sigma(w_i)|+|\sigma(w_{i+1})| \le M$, then
the closed path $\tau_i=\sigma(w_i) x_{i+1} \sigma(w_{i+1})^{-1}$ in \cay\
is of length $\le M+2$ and therefore represents up to cyclic reduction an
element of $R$.

     If $|\sigma (w_i)|+|\sigma (w_{i+1})| > M$ we break up the closed path
$\tau_i$ again, using the bicombing paths $\sigma_{i,t}=\sigma(\sigma(w_i)(t),
\sigma(w_{i+1})(t))$ that connect $\sigma(w_i)(t)$ to $\sigma(w_{i+1})(t)$
for all positive integers $t\le \max(|\sigma(w_i)|,|\sigma(w_{i+1})|)$.
By (\ref{konsk}),
\be\label{konsst}|\sigma_{i,t}|\le \max(2f(n/2)/k,M/2).\ee
Let $\sigma(w_i)=\avec{a}{p},\quad
\sigma(w_{i+1})=\avec{b}{q}$, $a_j ,b_l \in X^{\pm 1}\cup \{ 1\}$. We examine
the length of the closed paths $\tau_{i,t} $ that are generated by the
connecting paths
$\sigma_{i,t}$: $\tau_{i,t}=\sigma_{i,t} b_{t+1} \sigma_{i,t+1}^{-1}
a_{t+1}^{-1}$ (see fig.~\ref{abb1}).
\begin{figure}[ht]
\vspace{4cm}
\caption{A diagram for $w$} \label{abb1}
\end{figure}
If $|\sigma_{i,t}|+|\sigma_{i,t+1}|\le M$, then $|\tau_{i,t}|<M+2$ and
$\tau_{i,t}$ represents up to cyclic reduction an element in $R$.
Otherwise, we break $\tau_{i,t}$ up again using
bicombing paths $\sigma_{i,t,s}=\sigma(\sigma_{i,t}(s), \sigma_{i,t+1}(s))$
for $s\le \max(|\sigma_{i,t}|,|\sigma_{i,t+1}|)$.

   There is one exception, namely if we are close to the boundary. This
is because the path of length one between $w_i$ and $w_{i+1}$ is not
(necessarily) a combing line. But the condition $|\sigma (g,g)|\le Mk/2$
implies
$$|\sigma(\sigma(w_i,w_i)(0), \sigma(w_{i+1},w_{i+1})(0))|\le \max
\left( \frac{|\sigma(w_i,w_i)| + |\sigma(w_{i+1},w_{i+1})|}{k}, \frac{M}{2} \right)
\le M,$$
and we have a representation of an element of $R$.

   By (\ref{konsst}), $|\sigma_{i,t,s}|\le\max (4f(n/2)/k^2,M/2)$.
If $|\sigma_{i,t,s}|+|\sigma_{i,t,s+1}|\le M$
then the closed path $\tau_{i,t,s}$, using
$\sigma_{i,t,s}$, $\sigma_{i,t,s+1}^{-1}$
and the segments of length $\le 1$ along $\sigma_{i,t}$ and $\sigma_{i,t+1}$,
is of length $\le M+2$ and therefore represents an element in $R$.
Otherwise, we break up further in the same manner using connecting bicombing
paths of length $\le\max (8f(n/2)/k^3,M/2)$, etc.~until $2^df(n/2)/k^d \le M/2$.
In this way we find a van Kampen diagram for $w$. This proves that $G$ is
finitely presented. $d$ can be estimated as
the smallest integer greater or equal than $\log_{k/2} (2f(n/2)/M)$.

    The isoperimetric inequality has the form:
$$\mbox{\# (2--cells)}\le n\cdot (f(n/2)+1) \cdot 2(f(n/2)+2)/k \cdots
2^{d-1}(f(n/2)+2)/k^{d-1}\le$$
$$\frac{n (f(n/2)+2)^d 2^{d(d-1)/2}}{k^{d(d-1)/2}}=
n^{O(\log n)}$$
where $d$ is given as above.    \bewende

\noindent {\bf Remark:} 1. Condition (\ref{konsm}) is not necessary
in order to prove that the presentation is finite.\\
2. The growth of the isoperimetric function is faster than polynomial
but slower than exponential; therefore we call it sub--exponential.\\

\satzanf Each group that has a bicombing in the sense of \cite{Sh1} has a
bicombing of narrow shape. \satzende

\bewanf  By using the notation of the proof above, the bicombing in the
sense of Short is a narrow bicombing with $|\sigma_{i,t}|\le M/2$ and
$f(n)=mn$ for a given constant $m\in \maN$ and $d=1$ in this case. \bewende

\satzanf Let $P=\prae{X}{R}$ be a finite presentation for the group $G$ with a
bicombing of narrow shape $\sigma $ and let $f$ be the polynomial
from (\ref{konsm}) bounding $|\sigma (g)|$.\\
1. There is a polynomial isodiametric function for $P$ of the same degree
as $f$.\\
2. If $\sigma $ is geodesic, then the isodiametric function is linear.
\satzende

\bewanf  Let $w\in F$ be a reduced nontrivial word of length $n$, which is
trivial in $G$, and let $D$ be the van Kampen diagram
for $w$ constructed
in the proof of theorem \ref{satzIsop}. One can reach every vertex
in the diagram $D$ from the basepoint 1 by traveling part of a
bicombing line $\sigma (w_i)$ of the first generation then traveling
part of a bicombing line $\sigma_{i,t}$ of the second generation then
part of a bicombing line $\sigma_{i,t,s}$ of the third generation
etc.. The length of a bicombing line of the $l$-th generation is $\le 2^l
f(n/2)/k^l$, and the sum of the lengths of successive generations
 of bicombing
lines therefore is $\le f(n/2) (1+ 2/k + (2/k)^2 + ...) =  f(n/2)k/(k-2)$.
Hence $\frac{k}{k-2} f(n/2)$ is an isodiametric function for the presentation $P$.
If $\sigma $ is
geodesic, then $f$ is the identity and the above function is linear. \bewende

The next theorem follows an idea of M.~Bridson \cite{Bri}. It shows that
the definition of a bicombing of narrow shape cannot
be sharpened.

\satzanf Let $X$ be a finite generating set of the group $G$.
Choose for every pair $g,h\in G$ a geodesic $\sigma(g,h)\in$\cay. Then
$$\forall x,y\in X^{\pm 1},\qua\forall g,h\in G,\quad
|\sigma(\sigma(g,h)(t),\sigma(gx,hy)(t))|\le (|\sigma(g,h)|+|\sigma(gx,hy)|)/2+1 $$
holds for all integers $t\in [0,\infty [$. \satzende

\bewanf Let $C=(|\sigma(g,h)|+|\sigma(gx,hy)|)/2$. If $t\le C/2$, then following
$\sigma(g,h)$ backwards from $\sigma(g,h)(t)$ to $g$ then one edge to $gx$ and then going
to $\sigma(gx,hy)(t)$ along $\sigma(gx,hy)$ gives a path of length at most $C+1$.
For $t>C/2$ follow $\sigma(g,h)$ from $\sigma(g,h)(t)$ to the vertex $h$, then
go one edge to $hy$ and then to $\sigma(gx,hy)(t)$ backwards along $\sigma(gx,hy)$. This
gives a path of length at most $C$.   \bewende

\section{A class of Examples}

Let $P_q=\prae{x,y,z}{[x,y^q]=z, [x,z]=[y,z]=1}$ be a presentation
of the group $G_q$ where $q\ge 1$ and $[a,b]$ denotes the commutator
of $a$ and $b$. $G_1$ is the 3--dimensional integral
Heisenberg group. Let $F$ be the free group on $\{ x,y,z\} $. Let
$w,v\in F$. If both words are equal in $F$,
we write $w\equiv v$. If they are the same in $G_q$, we write $w=v$.

   It is easy to see, that
\be\label{kons1} z^{jl}=x^j y^{ql} x^{-j} y^{-ql} \ee
holds in $G_q$.

\lemanf[normal form for $G_q$] Let $w\in F$. Then, for $q>1$, there is a word

\be\label{nf} \tau (w)\equiv y^s x^{r_1} y^{s_1} x^{r_2} \ldots
y^{s_{m-1}} x^{r_m} y^p z^n \in F \ee
with $r_i,s_i\ne 0$ and
$$\mbox{for }q \quad \mbox{even:}\qquad s,s_i \in \{-q/2+1,\ldots ,q/2\},$$
$$\mbox{for }q \quad \mbox{odd:}\qquad s,s_i \in \{-(q-1)/2, \ldots , (q-1)/2 \} , $$
and, for $q=1$, there is a word
\be\label{nf2} \tau (w)\equiv x^r y^p z^n \in F \ee
such that $\tau (w)=w$ in $G_q$ and for all $v\in F$ with $w=v$ in $G_q$,
$\tau(w)\equiv \tau(v)$.    \lemende

\bewanf The case $q=1$ is trivial. For $q>1$ it is easy to see that each
word $w\in F$ can be transformed into $\tau (w)$ using the relations of $P_q$.
In order to prove uniqueness, let $w$ and $v$ be two words in $F$
representing the same element in $G_q$.
Let $H_q=G_q/ \normclz $, where
$\normclz $ denotes the normal closure of $z$
in $G_q$. $T_q=\prae{x,y}{xy^qx^{-1}=y^q}$ is a presentation for $H_q$,
which is an HNN--extension. Therefore
$w$ and $v$ have the same normal form (see \cite{LS}) $\tau'(w)=\tau'(v)$
in $H_q$ which is equal to the normal form in $G_q$,
except that $n=0$. Since $z$ is central, $\tau (w)$ and $\tau (v)$ can only
differ by a power of $z$. But $z$ has infinite order in $G_q$
which implies $\tau(v)\equiv \tau(w)$. \bewende

     The normal forms (\ref{nf}) and (\ref{nf2}) define a path
$\sigma (w)$ from 1 to $w$ in the Cayley
graph \cayq\ of $G_q$ for every $w\in F$. Define paths $\sigma (g,h)$
by taking equivariant lines; define

\be \sigma(g,h)(t):=g\cdot \sigma(1,g^{-1}h)(t) =g\cdot
\sigma(g^{-1}h)(t)\qquad         \forall g,h\in G_q \ee

\satzanf The paths $\sigma(g)$ are recursive
(i.e.~$|\sigma (g)|\le f(d(1,g))$) with a function
$f(x)=2x^2+3x$ for $q>1$ and $f(x)=x^2+x$ for $q=1$.    \satzende

\bewanf The relations in $P_q$ say that $z$ commutes with $x$ and $y$,
in particular any power of $z$ can be shifted to any place in a given word,
and that $x$ commutes with $y^q$ at the expense of introducing $z$ or $z^{-1}$.\\

For $q>1$, let $w \equiv \sigma(g)\equiv y^s x^{r_1} y^{s_1} x^{r_2} \ldots
y^{s_{m-1}} x^{r_m} y^p z^n \in F$ be the normal form for $g$.
We observe first that
\be\label{num1} d(1,g) \ge \sum_{i=1}^m |r_i|+ \sum_{i=1}^{m-1} |s_i| +|s|\ee
This is due to the fact that the exponents of the $y$-powers which occur
in $w$ can only be changed by adding multiples of $q$ (The relations
(\ref{kons1}) allow to permute powers of $x$ with powers of $y^q$). However,
the range for $s_i$ and $s$ in the normal form $w$ is such that $|s_i|$
and $|s|$ can not decrease under these changes. The same
argument also shows that\\
$d(1,g) \ge \sum_{i=1}^m |r_i|+ \sum_{i=1}^{m-1}
|s_i|+|s|+\max \{ |p|-(\sum |s_i| + |s|),0\} $, which implies:
\be\label{num2} d(1,g) \ge |p|\ee
Therefore $\sum_{i=1}^m |r_i|+ \sum_{i=1}^{m-1} |s_i| +|s|+|p|\le 2 d(1,g)$.
In order to prove $|w|=\sum_{i=1}^m |r_i|+ \sum_{i=1}^{m-1} |s_i| +|s|+|p|
+|n|\le f(d(1,g)) = 3 d(1,g) + 2 d^2(1,g)$,
we only need to show that  $|n|\le d(1,g) + 2 d^2(1,g)$:

        We claim that
\be\label{num3} d(1,g) \ge\sum_{i=1}^m |r_i|+ \sum_{i=1}^{m-1} |s_i| +|s|+\ee
$$\min \left\{ \max \left[ |n| - (\sum_{i=1}^m |r_i|+|r|)[(\sum_{i=1}^{m-1} |s_i|
+|s|+|p|)/q+|l|] ,0\right] + 2|r| + 2q|l| \right\} $$
where the minimum ranges over $|r|$ and $|l|$.
If $|n| \le (\sum_{i=1}^m |r_i|) (\sum_{i=1}^{m-1} |s_i|
+|s|+|p|)/q$ the minimum term on the right hand will be 0 and
the inequality holds by (\ref{num1}). If\\
$|n| > (\sum_{i=1}^m |r_i|) (\sum_{i=1}^{m-1} |s_i|+|s|+|p|)/q$
we observe first that $|n|$ may decrease by at most $|k||l|$ if a
power $y^{ql}$ is pushed across a power $x^k$ in $w$.

    If we do not introduce new powers of $x$ or $y^q$ by inserting $x^r x^{-r}$ or
$y^{ql} y^{-ql}$ into the word, the amount by which $|n|$ may be decreased by
means of permuting powers of $x$ with powers of $y^q$ is clearly bounded by
$\sum_{i=1}^m |r_i| (\sum_{i=1}^{m-1} |s_i|+|s|+|p|)/q$. This coarse estimate
stems from the following fact: Among all words in $x$ and $y$ whose sum of
absolute values of $x$-exponents and sum of absolute values of $y$-exponents
is the same as for $w$, $y^{\sum |s_i|+|s|+|p|} x^{\sum |r_i|}$ can absorb
the largest powers $z^{n'}$ or $z^{-n'}$ by permuting powers of $x$ with
powers of $y^q$.

      If we prolong the word by inserting $x^r x^{-r}$ and
$y^{ql} y^{-ql}$ at suitable places, the amount by which $|n|$ can be decreased
by means of (\ref{kons1}) is bounded by
$(\sum_{i=1}^m |r_i|+|r|)[ (\sum_{i=1}^{m-1} |s_i|+|s|+|p|)/q+|l|]$;
and, at the same time, the length of the x-y-part of the word increases
by $2|r|+2q|l|$. This explains inequality (\ref{num3}).

    Now, let $|r_0|$ and $|l_0|$ be the values for $|r|$ and $|l|$ for which the
minimum occurs in (\ref{num3}). Then
$d^2(1,g) \ge (\sum_{i=1}^m |r_i|+2|r_0|) (\sum_{i=1}^{m-1} |s_i|
+|s|+2|l_0|)$, and, by (\ref{num2}),
$d^2(1,g) \ge (\sum |r_i|+2|r_0|) |p|$
which implies $2 d^2(1,g) \ge (\sum_{i=1}^m |r_i|+|r_0|)[ (\sum_{i=1}^{m-1} |s_i|
+|s|+|p|)/q+|l_0|]$.
Therefore, by (\ref{num3}) again, $|n|\le d(1,g) + 2 d^2(1,g)$
which proves the Theorem for $q>1$.

    For $q=1$ the proof is similar, but easier. Let
$\sigma (g)\equiv x^r y^s z^n$. It is clear that $d(1,g)\ge |r|+|s|$.
If $|n|\le |r|+|s|$, then $d(1,g) + d(1,g)^2\ge |\sigma(g)|$;
if $|n| > |r|+|s|$, then, by the same ideas as in the proof for $q>1$,\\
$d(1,g) \ge |r|+|s|+ \min\{\max [|n|-(|r|+|r'|)(|s|+|s'|),0] +2|r'|+
2|s'|\}$ where the minimum ranges over the values of $|r'|$ and
$|s'|$. Let $|r_0'|$ and $|s_0'|$ be the values for which the minimum
occurs, then $|r|+|s|+|n| \ge d(1,g) + (|r|+|r_0'|)(|s|+|s_0'|)\ge
d(1,g) + d^2(1,g)$.     \bewende

\satzanf  $\sigma (g,h)$ defines a bicombing of narrow shape with constants \\
          $M=24q+18$ and $k=11/5$.  \satzende

\bewanf   Recall that a recursive $\sigma $ is of narrow shape, if
there exists an integer $M>1$ and a real $k>2$, such that
for all $g,h\in G$ and $a,b\in X^{\pm 1}\cup \{ 1\} $
$$ |\sigma (\sigma (g,h)(t),\sigma (ga,hb)(t))|
\le \max((|\sigma (g,h)|+|\sigma (ga,hb)|)/k,M/2) $$
holds for all integers $t\in [0,\infty [$. Since the bicombing is
equivariant, it suffices to show this inequality for $g=1$.

  For $q>1$ let $v\in F$ be in normal form $v\equiv y^s x^{r_1} y^{s_1} x^{r_2} \ldots
y^{s_{m-1}} x^{r_m} y^p z^n$, such that $v=h$ in $G_q$ ($\sigma(1,h)\equiv v$).
Let $w$ be the group element $a^{-1}vb$ brought into normal form
($\sigma(a,vb)\equiv w$) (see fig.~\ref{abb2}).
\begin{figure}[ht]
\vspace{4.5cm}
\caption{Close bicombing lines} \label{abb2}
\end{figure}

Now calculate the length of the bicombing lines (the {\it combing distance\/})
between these two paths $w,v$ in \cayq . Call the maximal combing distance
between two such paths \dist .

   If $a=1$ and $b\in \{ 1,z^{\pm 1}\} $, then \dist $\le 1$. If $a=1$
and $b\in \{ y^{\pm 1}\} $ then \dist $=2$.

    If $a=1$ and $b\in \{ x^{\epsilon }\} $ ($\epsilon = \pm 1$),
then \dist $\le |l|+q+1$, where $l$ is such that $-q/2+1\le p-lq\le q/2$ for
$q$ even and $-(q-1)/2\le p-lq\le (q-1)/2$ otherwise.
To see this, observe that $v$ ends with $y^p z^n$ but $w$ ends
with $x^\epsilon y^{ql} z^{n-\epsilon l}$. Since $|w|+|v|\ge 2q|l|$ we get
for $q\ge 2$
and \dist $> M/2$: $(|w|+|v|)/k> $\dist .

  There are a few more cases which are relatively easy. The most critical case
which requires the sharpest estimates occurs if
$a=y^\epsilon ,\qua b=x^\alpha $ with $\alpha ,\epsilon\in\{ \pm 1\} $; in
particular if $y^s$ is at the boundary of its range to which it is restricted
by the normal form, and the premultiplication by $a^{-1}=y^{-\epsilon} $ moves
it out of this range, as, for example, in the
case $\epsilon =-1,\qua s=q/2$ and $q$ even (the
other cases can be treated similarly).

    In this case
$v\equiv y^{q/2} x^{r_1} y^{s_1} x^{r_2} \ldots y^{s_{m-1}} x^{r_m} y^p z^n$
and $$w\equiv y^{-q/2+1} x^{r_1} y^{s_1} x^{r_2} \ldots y^{s_{m-1}} x^{r_m}
y^{p-lq} x^\alpha y^{(l+1)q} z^{n-\sum r_i-\alpha (l+1)},$$
where $l$ is as above. Using the rule $|a|+|a-b|\ge |b|$ we obtain the estimate:
$|w|+|v|\ge 2q|l|+2\sum |r_i|+|\sum r_i+\alpha (l+1)|$.
A careful study of the lengths of the combing distances shows that
$$\mbox{\dist} \le \max\{ \max_{j\le m} |\sum_{i=1}^j r_i|+1,
|\sum_{i=1}^m r_i|+|l|+3q+2\}\le (\sum_{i=1}^m |r_i|+|\sum_{i=1}^m r_i|)/2
+|l|+3q+2.$$
Since $q\ge2$ and $k=11/5$,
$(|w|+|v|)/k\ge 20|l|/11+10\sum |r_i|/11+5|\sum r_i+\alpha (l+1)|/11$.
We will show that the right hand side is
$\ge (\sum_{i=1}^m |r_i|+|\sum_{i=1}^m r_i|)/2+|l|+3q+2$
whenever \dist $>M/2$ (which, by the above estimate for \dist ,
proves the Theorem for this case). This is equivalent to:
$$9|l|+10\sum |r_i| +5|\sum r_i+\alpha (l+1)|\ge 11\sum |r_i|/2+11|\sum r_i|/2+33q+22.$$
The left hand side can be simplified by the following estimates:
$5|l|+5|\sum r_i +\alpha +\alpha l|\ge 5|\sum r_i +\alpha |\ge 5|\sum r_i|-5$,
and $10\sum |r_i|+5|\sum r_i|\ge 19\sum|r_i|/2+11|\sum r_i|/2\ge 4\sum |r_i|+11|\sum r_i|/2$.
Therefore the above inequality follows from $4(|l|+\sum |r_i|)\ge 33q+27$,
which follows from \dist $> M/2$ using the value $M=24q+18$ and the estimate
\dist $\le (\sum |r_i|+|\sum r_i|)/2+|l|+3q+2\le \sum |r_i|+3q+2$.

    The proof for $q=1$ is much simpler and left to the reader.    \bewende

In the following we use Cockcroft 2-complexes to get lower bounds for
isoperimetric functions. This idea is due to S.~Gersten \cite{G9}.

\satzanf $G_q$ has no quadratic isoperimetric inequality and therefore no
combing in the sense of Short \cite{Sh1}. \satzende

\bewanf  There is a van Kampen diagram for $w_n\equiv [x^n,y^{qn}]\cdot [y^{-qn},x^{-n}]$
in $G_q$, which
has $n^3$ more 2-cells $[x,z]$ of positive then of negative type.
W.~A.~Bogley proved
in \cite{B1}, that the corresponding 2-complex is Cockcroft. So each
$\pi_2$--element has the same number of positive as
of negative 2-cells $[x,z]$, which proves that every van Kampen diagram for
$w_n$ will contain at least $n^3$ 2-cells $[x,z]$ and so proves
the theorem.    \bewende

\bibliographystyle{plain}

\vspace{2cm}

\begin{tabbing}
\hspace*{1.5cm} \= Stephan Rosebrock \hspace{5cm} \= G\"unther Huck\\
\> Institut f.~Didaktik der Mathematik \> Dept. of Math.\\
\> J.- W.- Goethe Universit\"at           \> Northern Arizona University\\
\> Senckenberganlage 9                  \> Flagstaff AZ 86011\\
\> 60054 Frankfurt/M\@.                  \> USA\\
\> West-Germany
\end{tabbing}

e-mail:\\
huck@nauvax.ucc.nau.edu\\
rosebrock@mathematik.uni-frankfurt.d400.de

\end{document}